\documentclass[11pt]{article}
\usepackage{amsmath}
\usepackage{amsfonts}
\usepackage{amssymb}
\usepackage{graphicx}

\begin{document}

\newcommand{\ls}[1]
   {\dimen0=\fontdimen6\the\font \lineskip=#1\dimen0
\advance\lineskip.5\fontdimen5\the\font \advance\lineskip-\dimen0
\lineskiplimit=.9\lineskip \baselineskip=\lineskip
\advance\baselineskip\dimen0 \normallineskip\lineskip
\normallineskiplimit\lineskiplimit \normalbaselineskip\baselineskip
\ignorespaces }
\renewcommand{\theequation}{\thesection.\arabic{equation}}

\newtheorem{definition}[equation]{Definition}
\newtheorem{lemma}[equation]{Lemma}
\newtheorem{proposition}[equation]{Proposition}
\newtheorem{corollary}[equation]{Corollary}
\newtheorem{example}{Example}[section]
\newtheorem{conjecture}{Conjecture}[section]
\newtheorem{algorithm}{Algorithm}[section]
\newtheorem{theorem}[equation]{Theorem}
\newtheorem{exercise}[equation]{Exercise}
\newtheorem{assumption}[equation]{Assumption}

\def\le{\leq}
\def\ge{\geq}
\def\lt{<}
\def\gt{>}

\newcommand{\Exercise}{{\bf Exercise}}
\newcommand{\req}[1]{(\ref{#1})}
\newcommand{\beps}{{\scriptscriptstyle{{ \cal E}}}}
\newcommand{\lip}{\langle}
\newcommand{\rip}{\rangle}
\newcommand{\lf}{\lfloor}
\newcommand{\lc}{\lceil}
\newcommand{\rc}{\rceil}
\newcommand{\rf}{\rfloor}
\newcommand{\supp}{{\rm supp\,}}
\newcommand{\ARn}
 {\begin{array}[t]{c}
  \longrightarrow \\[-0.3cm]
 \scriptstyle {n\rightarrow \infty}
  \end{array}}

\newcommand{\uu}{\underline}
\newcommand{\oo}{\overline}
\newcommand{\dfr}{\displaystyle\frac}
\newcommand{\La}{\Lambda}
\newcommand{\la}{\lambda}
\newcommand{\eps}{\varepsilon}
\newcommand{\om}{\omega}
\newcommand{\Om}{\Omega}
\newcommand{\Inv}{{\text{\rm Inv}\,}}
\newcommand{\crit}{{\rm crit}}

\newcommand{\EE}{{\mathbb E}}
\newcommand{\NN}{{\mathbb N}}
\newcommand{\PP}{{\mathbb P}}
\newcommand{\QQ}{{\mathbb Q}}
\newcommand{\reals}{{\mathbb R}}
\newcommand{\ZZ}{{\mathbb Z}}

\newcommand{\calA}{{\mathcal A}}
\newcommand{\calB}{{\mathcal B}}
\newcommand{\calF}{{\mathcal F}}
\newcommand{\calG}{{\mathcal G}}
\newcommand{\calI}{{\mathcal I}}
\newcommand{\calO}{{\mathcal O}}
\newcommand{\calQ}{{\mathcal Q}}
\newcommand{\calS}{{\mathcal S}}
\newcommand{\calV}{{\mathcal V}}
\newcommand{\calW}{{\mathcal W}}
\newcommand{\calY}{{\mathcal Y}}
\newcommand{\calX}{{\mathcal X}}

\newcommand{\bfcdot}{{\boldsymbol \cdot}}

\newcommand{\won}{{\boldsymbol 1}}
\newcommand{\hilbert}{\bigcirc\kern -0.8em
              {\rm\scriptstyle {H}\;}} 

\newcommand{\limn}{\lim_{n \rightarrow \infty}}
\newcommand{\limk}{\lim_{k \rightarrow \infty}}
\newcommand{\limi}{\lim_{i \rightarrow \infty}}
\newcommand{\liml}{\lim_{\ell \rightarrow \infty}}
\newcommand{\limv}{\lim_{v \rightarrow \infty}}
\newcommand{\limm}{\lim_{m \rightarrow \infty}}
\newcommand{\limd}{\lim_{\delta \rightarrow \infty}}
\newcommand{\limsupn}{\limsup_{n \rightarrow \infty}}
\newcommand{\liminfn}{\liminf_{n \rightarrow \infty}}

\newcommand{\proof}{\noindent \textbf{Proof:\ }}
\newcommand{\remark}{\noindent \textbf{Remark:\ }}
\newcommand{\remarks}{{\bf Remarks:}}

\newcommand{\eff}{{\operatorname{eff}}}

\newcommand{\dsum}{\displaystyle\sum}
\newcommand{\dprod}{\displaystyle\prod}

\newcommand{\ffrac}[2]
  {\left( \frac{#1}{#2} \right)}

\newcommand{\one}{\frac{1}{n}\:}
\newcommand{\half}{\frac{1}{2}\:}

\def\squarebox#1{\hbox to #1{\hfill\vbox to #1{\vfill}}}
\newcommand{\qed}{\hspace*{\fill}
           \vbox{\hrule\hbox{\vrule\squarebox{.667em}\vrule}\hrule}\smallskip}


\title{{Aging properties of Sinai's model of random walk in random environment}}
\author{Amir Dembo\\
Department of Statistics and\\
Department of Mathematics\\
Stanford University\\
Stanford, CA 94305, USA \and Alice Guionnet
\\
Ecole Normale Superieure de Lyon\\
Unite de Mathematiques pures et appliquees\\
UMR 5669, 46 Allee d'Italie \\
69364 Lyon Cedex 07, France
\and
Ofer Zeitouni\\
Department of Electrical Engineering\\
Technion\\
Haifa 32000, Israel}

\date{{\textit{May 26, 2001}}}
\maketitle


\ls{1}

\abstract{We study in this short note aging properties of  Sinai's
(nearest neighbour) random walk in random environment. 
With $\PP^o$ denoting the annealed law of the
RWRE $X_n$, our main result is a full proof of the
following statement due to 
P. Le Doussal, C.  Monthus and D. S.   Fisher:
$$\lim_{\eta\to0} \lim_{n\to\infty} \PP^o
\left(\frac{|X_{n^h} - X_n|}{(\log n)^2} < \eta\right)
= \frac{1}{h^2} \left[ \frac{5}{3} - \frac{2}{3} e^{-(h-1)} \right].
$$
}
\section{Introduction: Sinai's model and non standard limit laws}
\setcounter{equation}{0}
\label{sec-sinai}

We shall here study Sinai's random walk in random 
environment, and focus on its aging properties. We begin
by recalling Sinai's model and it's long time behaviour.
This review is not meant to be exhaustive: see the recent
paper \cite{shi} for an updated account of Sinai's model which includes
many topics not covered here.

Let
$N_z = \{\omega_z^-,\omega_z^0,\omega_z^+ \geq 0 :
\omega_z^- + \omega_z^0 +\omega_z^+ =1\}$ for $z \in \ZZ$,
and $\Omega = \prod_{z \in \ZZ} N_z$ equipped with the 
product topology and the corresponding Borel $\sigma$-field
and shift operator. The {\it random
environment} is $\omega \in \Omega$ of
law $P \in M_1(\Omega)$, a Borel probability measure on $\Omega$.
For each $\omega \in \Omega$ we define (Sinai's) {\it random walk in
random environment} as the time-homogeneous Markov
chain $\{ X_n \}$ taking values in $\ZZ$ with transition
probabilities
$$
P_\omega(X_{n+1}=z+1|X_n=z) = \omega_z^+, \quad
P_\omega(X_{n+1}=z|X_n=z) = \omega_z^0, \quad
P_\omega(X_{n+1}=z-1|X_n=z) = \omega_z^- \;.
$$
We use $P_\omega^v$ to denote the law
induced on $\ZZ^\NN$ when $P_\omega^v(X_0=v)=1$, 
refering to it as the {\it quenched} law of $\{X_n\}$.
Noting that $\omega \mapsto P_\omega^v(G)$ is Borel measurable 
for any fixed Borel set $G \subset \ZZ^\NN$, we
define the {\it annealed} law of $\{ X_n \}$ as
$P^v(G)=\int_\Omega P_\omega^v(G) dP(\omega)$. Note that
under $P^v$ the random walk in random environment $\{ X_n \}$
is not a Markov chain!

With $\rho_z=\omega_z^-/\omega_z^+$ and
$\oo{R}_k= k^{-1} \sum_{i=0}^{k-1} \log \rho_i$, we 
assume throughout that the following holds.
\begin{assumption}
\label{ass-sinai}
The probability measure $P$ is stationary,
strongly mixing on $\Omega$ and such that 
$P(\om_0^+ + \om_0^- > 0) = 1$, 
$E_P \log \rho_0 = 0$, and there exists an $\varepsilon > 0 $ such that
$E_P|\log \rho_0|^{2+\varepsilon} < \infty$. 
Further,
$\{\sqrt{k}\, \oo R_{[kt]}/\sigma_P\}_{t\in \reals}$ converges weakly to
a Brownian motion for some 
$\sigma_P> 0$.
\end{assumption}
(In the i.i.d. case, $\sigma_P^2=E_P (\log \rho_0)^2$).
Define 
$$
W^n (t) = \frac{1}{\log\,n} \sum_{i=0}^{\lfloor (\log n)^2 t \rfloor}
\log \rho_i \cdot (\text{sign\ } t)
$$
with $t\in\reals$.  By Assumption~\ref{ass-sinai}, 
$\{W^n(t)\}_{t\in\reals} $ converges weakly to $\{\sigma_P B_t\}$, where
$\{B_t\}$ is a two sided Brownian motion.

Next, we call a triple $(a,b,c)$ with $a < b < c$ a valley of the path
$\{W^n (\cdot) \}$ if
\begin{align*}
W^n (b) & = \min_{a \le t \le c} W^n (t)\,,\\
W^n (a) & = \max_{a \le t \le b} W^n (t)\,,\\
W^n (c) & = \max_{b \le t \le c} W^n (t)\,.
\end{align*}
The \textsl{depth} of the valley is defined as
$$
d_{(a,b,c)} = \min (W^n(a) - W^n(b), W^n(c) - W^n(b))
\,.
$$
If $(a,b,c)$ is a valley, and $a<d< e < b$
are such that
$$
W^n (e) - W^n(d) = \max_{a\le x < y\le b} W^n(y) - W^n(x)
$$
then $(a,d,e)$ and $(e,b,c)$ are again valleys, which are obtained from
$(a,b,c)$ by a \textsl{left refinement}.
One defines similarly a 
\textsl{right refinement}.
Define
\begin{align*}
c_0^n & = \min \{t\ge 0:\quad W^n(t) \ge 1\}\\
a_0^n & = \max \{t \le 0:\quad W^n (t) \ge 1\}\\
W^n(b_0^n) & = \min_{a_0^n\le t \le c_0^n} W^n(t)
\,.
\end{align*}
($b_0^n$ is not uniquely defined, however, due to Assumption \ref{ass-sinai},
with $P$-probability approaching 1 as $n\to\infty$,  
all candidates for $b_0^n$
are within distance converging to 0 as $n\to\infty$; we define $b_0^n$ then as
the smallest one in absolute value.)

\begin{figure}[h]
\begin{center}
\includegraphics[scale=.6]{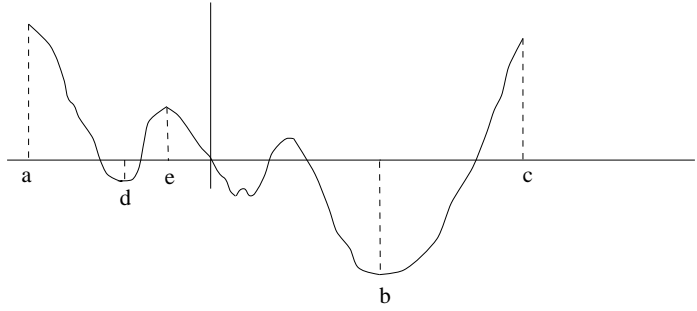}
\caption{Left refinement of $(a,b,c)$}
\end{center}
\end{figure}

One may now apply a (finite) sequence of refinements to find the
\textit{smallest} valley $(\oo{a}^n, \oo{b}^n, \oo{c}^n)$
with $\oo{a}^n < 0 < \oo{c}^n$, while
$d_{(\oo{a}^n, \oo{b}^n, \oo{c}^n)} \geq 1$. 
We define similarly the smallest valey 
$(\oo{a}^n_\delta, \oo{b}^n_\delta, \oo{c}^n_\delta)$ 
such that 
$d_{(\oo{a}^n_\delta, \oo{b}^n_\delta, \oo{c}^n_\delta)} \geq 1+\delta$.
Let
\begin{equation*}
A_n^{J,\delta} = 
\Biggl\{ \om\in \Omega: 
\begin{array}[t]{c}\oo b^n=\oo b^n_\delta, \,
\text{any refinement $(a,b,c)$ of $(\oo{a}^n_\delta, \oo{b}^n_\delta, 
\oo{c}^n_\delta)$ with $b\neq \oo{b}^n$ 
has depth } < 1-\delta\,,\\
\min_{t\in [\oo{a}^n, \oo{c}^n]\backslash [\oo{b}^n - \delta, \oo{b}^n +
\delta]}
W^n(t) - W^n (\oo b^n) > \delta^3\,,
|\oo{a}^n_\delta| + |\oo{c}^n_\delta| \le J
\end{array}
\Biggr\}
\end{equation*}
then it is easy to check by the properties of Brownian motion that
\begin{equation}
\label{sinai1}
\lim_{\delta \to 0}
\lim_{J\to\infty} 
\lim_{n\to\infty}
P(A_n^{J,\delta}) = 1\,.
\end{equation}
The following theorem is due to 
\cite{sinai}. For completeness, we include
a proof, which follows the approach 
of \cite{golosov}, who dealt with a RWRE reflected at $0$,
i.e. with state space $\ZZ_+$. 
\begin{theorem}
\label{the-sinai}
Assume $P(\min(\om_0, \om_0^+) < \eps) = 0$ and Assumption~\ref{ass-sinai}.
For any $\eta > 0$,
$$
\PP^o \left( \left| \frac{X_n}{(\log n)^2} - \oo{b}^n \right| > \eta
\right) \underset{n\to\infty}{\to}0
$$
\end{theorem}
\proof
Fix $\delta < \eta/2, J$ and $n_0$ large enough such that for all $n>n_0$,
$\om\in A_n^{J,\delta}$.
For simplicity of notations, assume in the sequel that
$\om$ is such that $\oo{b}^n > 0$.
Write $a^n=\oo a^n (\log n)^2, b^n=\oo b^n (\log n)^2,
c^n=\oo c^n (\log n)^2$, with similar notations for 
$a^n_\delta,b^n_\delta,c^n_\delta$. Define
\begin{align*}
\oo{T}_{b,n} & = \min \{ t \ge 0: 
X_t  = {b}^n \,
\text{or\ } X_t=  {a}^n_\delta \} 
\,.
\end{align*}

We next recall that in dimension one, harmonic functions are readily 
evaluated.
For $z \in [-m_-, m_+]$, define
$$
\calV_{m_-, m_+, \om} (z) :=
P_\om^z (\{X_n\} \text{\ hits\ $-m_-$ before hitting $m_+$)}
\,.
$$
The Markov property implies that 
$\calV_{m_-, m_+, \om} (\cdot)$ is harmonic, that is it satisfies
\begin{equation}
\label{1}
\left\{
\begin{aligned}
(\om_z^++\om_z^-)\calV_{m_-, m_+, \om}(z) & = \om_z^- \calV_{m_-,m_+, \om} (z-1)\\
& \quad +\om_z^+ \calV_{m_-, m_+, \om} (z+1), \quad
            z\in (-m_-, m_+) \\
\calV_{m_-,m_+, \om}(-m_-) & =1, \quad  \calV_{m_-,m_+, \om}(m_+) =0
\end{aligned}
\right.
\end{equation}
Solving \req{1} (noting that the solution is unique due to
the maximum principle), we find
\begin{equation}
\label{2}
\calV_{m_-, m_+, \om} (z) = 
\frac{\dsum_{i=z+1}^{m_+}\; 
\dprod_{j=z+1}^{i-1} \rho_j}
{\dsum_{i=z+1}^{m_+} \; \dprod_{j=z+1}^{i-1} \rho_j + \dsum_{i=-m_-+1}^z
\left(\dprod_{j=i}^z \rho_j^{-1}\right)}
\end{equation}
Hence, for $\om\in A_n^{J,\delta}$, 
\begin{equation}
\label{sinai2}
P_\om^o \Bigl(X_{\oo{T}_{b,n}} 
= \calV_{a_\delta^n, b^n_\delta, \om}(0)={a}^n_\delta \Bigr)
 \le \frac{1}{1+ \frac{\exp\{
(\log n)(W^n(\oo a^n_\delta) - W^n(\oo b^n))\}}{J(\log n)^2}}
 \le \frac{J(\log n)^2}{n^{1+\delta}}
\end{equation}
On the other hand, let $\tilde T_{b,n}$ have the law of 
$\oo{T}_{b,n}$ except that the walk $\{X_\cdot\}$ is reflected at
${a}^n_\delta $, and define similarly $\tilde \tau_1$.  
Using the Markov property, c.f. \cite{GZrev} for a similar computation,
we have that
$$
E_\om^o (\tilde 
\tau_1) = \frac{1}{\om_0^+} + \frac{\rho_0}{\om_{(-1)}^+} +
\cdots + \frac{ \prod_{i=0}^{ {a}^n_\delta + 2} \rho_{-i}}
{\om_{{a}^n_\delta-1}^+} + \prod_{i=0}^{{a}^n_\delta + 1} \rho_{-i}\,.
$$
Hence, with
$\tilde\om_i = \om_i$ for 
$i\not={a}^n_\delta$ and $\tilde\om_{{a}^n_\delta}^+ =
1$, for all $n$ large enough, 
\begin{align*}
E_\om^o(\oo{T}_{b,n}) & \le E_\om^o (\tilde{T}_{b,n})  =
\sum_{i=1}^{{b}^n} \sum_{j=0}^{i-1-{a}^n_\delta}
\frac{\prod_{k=1}^j \rho_{i-k}}{\om_{(i-j-1)}^+} \\
& \le \frac{1}{\eps} \sum_{i=1}^{{b}^n} \sum_{j=0}^{i-1-{a}^n_\delta}
e^{(\log n) (W^n(i) - W^n(i-j))} 
 \le \frac{2J^2}{\eps} e^{\log n (1-\delta)} \le
n^{1- \frac{\delta}{2}}\,.
\end{align*}
We thus conclude with (\ref{sinai2}) 
that
$$
	P_\om^o \Bigl( \oo{T}_{b,n} < n, \quad X_{\oo{T}_{b,n}} =
{b}^n\Bigr) \underset{n\to\infty}{\longrightarrow} 1
$$
implying that
\begin{equation}
\label{sinai3}
P_\om^o \Bigl( T_{{b}^n} < n \Bigr) \underset{n\to\infty}{\to} 1\,.
\end{equation}
Next note that another application of \req{2} yields
\begin{align}
\label{sinai4}
P_\om^{{b}^{n}{-1}} (X_\cdot\, \text{ 
hits ${b}^n$ before ${a}^n_\delta$})
& \ge 1- n^{-(1+\frac{\delta}{2})}\notag \\
P_\om^{{b}^{n}{+1}} (X_\cdot\, \text{
hits ${b}^n$ before ${c}^n_\delta$})
& \ge 1- n^{-(1+\frac{\delta}{2})} 
\end{align}
On the same probability space,
construct a RWRE $\{\tilde{X}_t\}$ with the same transition mechanism 
as
$\{X_t\}$ except that it is reflected at
${a}^n_\delta$, i.e.\ replace
$\om$ by $\tilde \om$.
Then, using (\ref{sinai4}),
\begin{align*}
P_\om^o  \left( \left| \frac{X_n}{(\log n)^2} -
\oo{b}^n \right| > \delta
\right)
& \le P_\om^o \Bigl(T_{{b}^n} > n \Bigr)
+ \max_{t\le n} P_\om^{{b}^n} \left(\left| \frac{X_t}{(\log n)^2} -
\oo{b}^n \right| > \delta \right) \\
&\le P_\om^o \Bigl(T_{{b}^n} > n\Bigr) +  \Bigl[
1-(1-n^{-(1+\frac{\delta}{2})})^n\Bigr]
+ \max_{t\le n} P_\om^{{b}^n} \left(\left| \frac{\tilde{X}_t}{(\log n)^2}
- \oo{b}^n \right| > \delta\right)
\end{align*}
Hence, in view of \req{sinai1} and \req{sinai3}, 
the theorem holds as soon as
we show that for $\om$ as considered here,
\begin{equation}
\label{sinai5}
 \max_{t\le n} P_\om^{\oo{b}^n} \left(\left| \frac{\tilde{X}_t}{(\log n)^2}
- \oo{b}^n \right| > \delta\right)
\underset{n \to \infty}{\longrightarrow} 0
\end{equation}
To see (\ref{sinai5}), define
$$
f(z) = \frac{\prod_{{a}^n_\delta +1 \le i < z} \om_i^+}
{\prod_{{a}^n_\delta +1 \le i < z} \om_{i+1}^-} ,\quad
\oo{f}(z) = \frac{f(z)}{f({b}^n)}
$$
(as usual, the product over an empty set of indices is taken as 1.
$\oo{f}(\cdot)$ 
corresponds to the invariant measure for the resistor network
corresponding to $\tilde{X}_.$).  Next, define the operator
\begin{equation}
\label{sinai6}
(Ag) (z)  = \oo{\om}_{z-1}^{\,+} g(z-1) + \oo\om_{z+1}^{\,-}
g(z+1) + \oo\om_z^{\,0}g(z)
\end{equation}
where $\oo{\om}_z = \om_z$ for
$z> {a}^n_\delta$,
$\oo{\om}_{{a}^n_\delta}^+ = 1, \oo\om_{{a}^n_\delta-1}^+ = 0$.
Note that
$A\oo{f} = \oo{f}$, and further that
$$
P_\om^{{b}^n} (\tilde X_t = z) = A^t {\bf 1}_{b^n}(z)\,.
$$
Since $\oo{f}(z) \ge {\bf 1}_{b^n}(z)$ and $A$ is a positive operator,
we conclude that
$$P_\om^{{b}^n} (\tilde{X}_t = z) \le \oo{f}(z)\,.
$$
But, for $z$ with $|z/(\log n)^2 - \oo{b}^n | > \delta$
and $\om\in A^{J,\delta}_n$, it holds that
$\oo{f} (z) \le e^{-\delta^3 \log n}$, and hence
$$
P_\om^{{b}^n} (\tilde{X}_t=z) \le n^{-\delta^3}\,.
$$
Thus,
$$
\max_{t\le n} P_\om^{\oo{b}^n} 
\left(\left| \frac{\tilde{X}_t}{(\log n)^2} - \oo{b}^n
\right| > \delta \right) \le n^{-\delta^3}\,,
$$
yielding (\ref{sinai5}) and completing the proof of the theorem.
\qed

We next turn to a somewhat more detailed study of the random variable
$\oo{b}^n$.  By replacing $1$ with $t$ in the definition of $\oo{b}^n$, one
obtains a process
$\{\oo{b}^n(t)\}_{t\ge 0}$.  Further, due to Assumption~\ref{ass-sinai},
the process $\{\oo{b}^n (t/\sigma_P)\}_{t\ge 0}$ 
converges weakly to a process
$\{\oo{b}(t)\}_{t\ge 0}$, defined in terms of the Brownian motion
$\{B_t\}_{t\ge 0}$; Indeed, $\oo{b}(t)$ is the location of the bottom of the
smallest valley of
$\{B_t\}_{t\ge 0}$, which surrounds $0$ and has depth $t$. Throughout
this section we denote by $\calQ$ the law of the 
Brownian motion $B_\cdot$. Our next goal is
to characterize the process
$\{\oo{b}(t)\}_{t\ge0}$.  Toward this end,
define
\begin{align*}
m_+ (t) &= \min\{B_s: 0 \le s \le t\}\,,
\quad
m_- (t) = \min\{B_{-s}: 0 \le s \le t\} \\
T_+ (a) & = \inf \{s\ge 0: B_s-m_+(s) = a\}\,,
\quad
T_- (a)  = \inf \{s\ge 0: B_{-s}-m_-(s) = a\}\\
s_\pm (a) & = \inf\{s\ge 0: m_\pm (T_\pm(a)) = B_{\pm s} \} \,,
\quad
M_\pm (a)  = 
\sup \{ B_{\pm \eta}: \quad 0\le \eta \le s_\pm (a) \} \,.
\end{align*}
Next, define $W_\pm(a) = B_{s_\pm (a)} $.
It is not hard to check that the pairs
$(M_+(\cdot), W_+ (\cdot))$ and $(M_-(\cdot), W_-(\cdot))$ form  independent
Markov processes.
Define finally
$$
H_\pm(a) = (W_\pm (a) + a) \vee M_\pm (a)\,.
$$

\begin{figure}[h]
\begin{center}
\includegraphics[scale=.6]{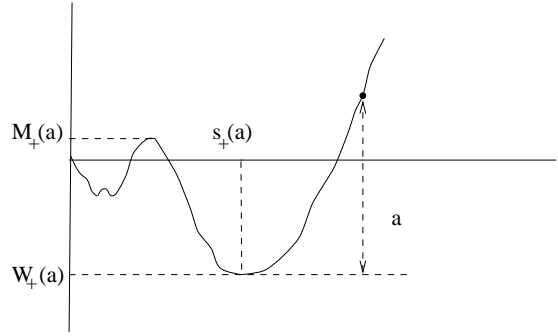}
\caption{The random variables $(M_+(a), W_+ (a), s_+(a))$}
\end{center}
\end{figure}

The following characterization of the law of $\oo{b}(a)$ 
is due to  \cite{kesten1}. This characterization can be found
also in \cite{golosov2}.
\begin{theorem}
\label{the-sinai2}
For each $a>0$, $\calQ(\oo{b}(a) \in \{s_+(a), -s_-(a)\})=1$.
Further, $\oo{b}(a) = s_+(a)$ iff
$H_+(a) < H_- (a)$.
\end{theorem}

\proof
Note that $\calQ(H_+ (a) = H_-(a))=0$.
That $\oo{b}(a) \in \{s_+(a), -s_-(a)\}$ is a direct consequence of the
definitions, i.e.\ assuming $\oo{b}(a) > 0$ and $\oo{b} (a) 
\neq s_+(a)$ it is easy to show
that one may refine from the right
the valley defining $\oo{b}(a)$, contradicting 
minimality.  We begin by showing, after Kesten~\cite{kesten1}, that
$\oo{b}(a) = s_+(a)$ iff either
\begin{equation}
\label{sinai7}
W_-(a) > W_+ (a), \quad M_+(a) < (W_-(a) + a) \vee M_-(a)
\end{equation}
or
\begin{equation}
\label{sinai8}
W_-(a) < W_+ (a), \quad M_-(a) > (W_+(a) + a) \vee M_+(a)
\,.
\end{equation}
Indeed, assume $\oo{b} (a) = s_+(a)$, and
$W_-(a) > W_+ (a)$.
Let $(\alpha, \oo{b}(a),\gamma)$ denote the minimal valley defining 
$\oo{b}(a)$.
If $-s_-(a) \le \alpha$, then
\begin{equation}
\label{sinai9}
M_-(a) = \max \{B_{-s}: \; s\in (0,s_-(a))\} \ge B_{-\alpha}
= \max \{B_{s}: -\alpha \le s \le \oo{b} (a) \} \ge M_+(a)
\end{equation}
implying (\ref{sinai7}).  On the other hand, if
$-s_-(a) > \alpha$, refine $(\alpha, \oo{b} (a), \gamma)$ on the left
(find $\alpha', \beta'$ with $\alpha<\alpha' < \beta' < \oo{b} (a)$),
such that
$$
B_{\beta'} - B_{\alpha'} = \max_{\alpha < x< y < \oo{b} (\alpha)}
(B_y - B_x) \ge M_+(a) - W_-(a)
$$
and thus minimality of $(\alpha, \oo{b}(a), \gamma)$ implies that
$M_+(a) - W_-(a) < a$, implying (\ref{sinai7}).

\begin{figure}[h]
\begin{center}
\includegraphics[scale=.5]{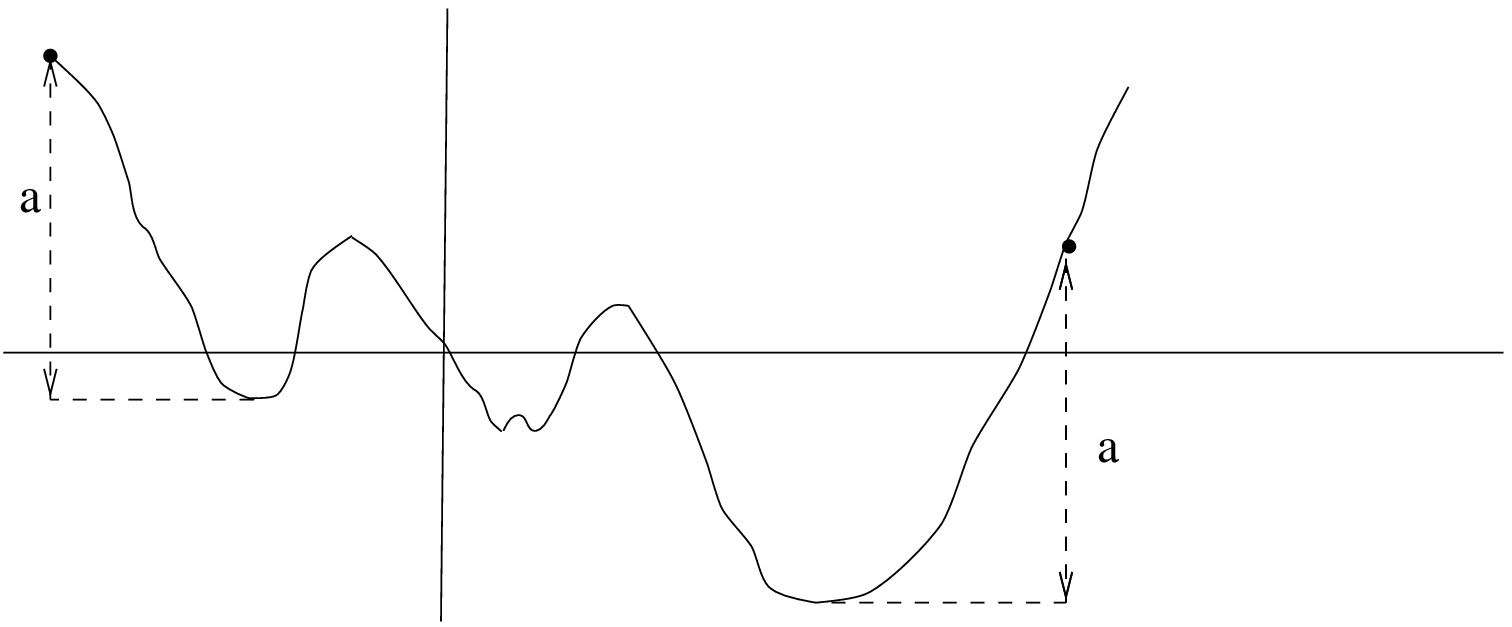}
\qquad \quad
\includegraphics[scale=.5]{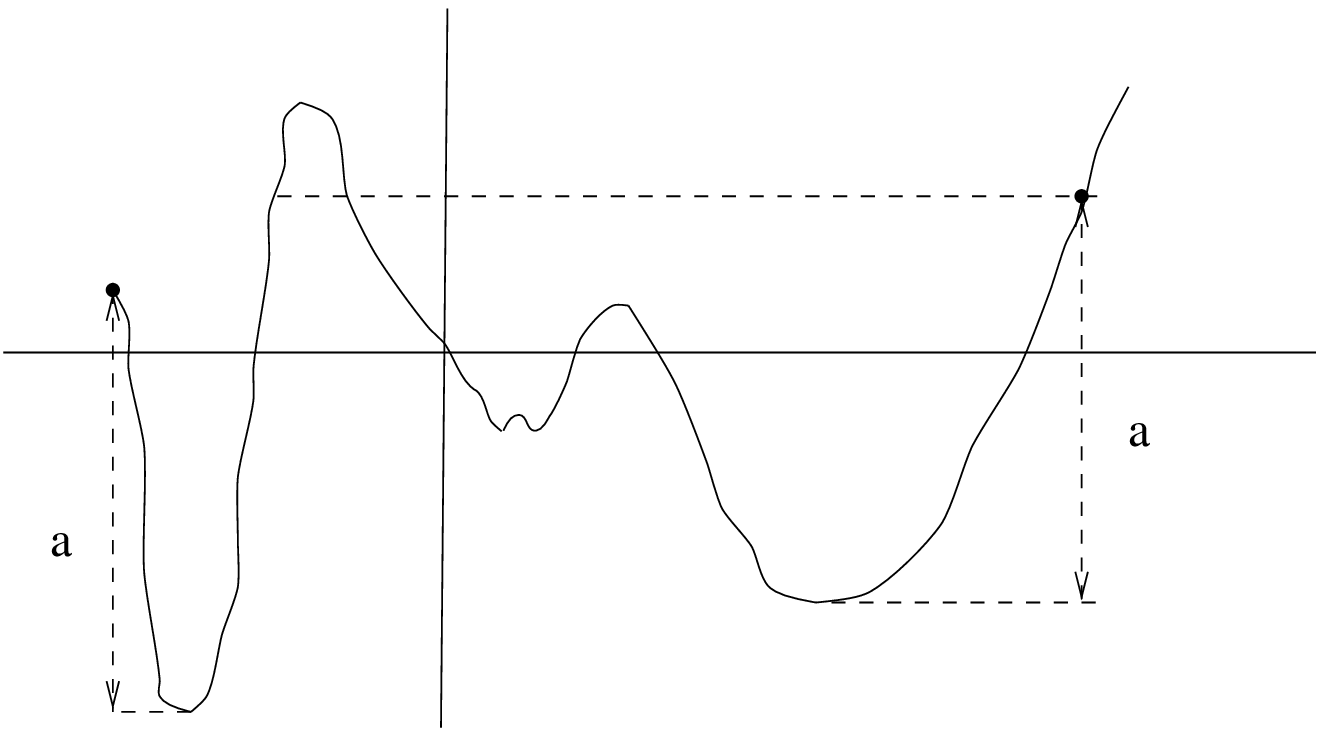}
\caption{$\oo{b}(a) = s_+(a)$}
\end{center}
\end{figure}

We thus showed that if $\oo{b}(a) = s_+(a)$ and
$W_-(a) > W_+(a)$ then (\ref{sinai7}) holds.  On the other hand, if
(\ref{sinai7}) holds, we show that $\oo{b}(a)=s_+(a)$ by considering the cases
$\alpha \le -s_-(a)$ and $-s_-(a) < \alpha$ separately. 
In the former case,
necessarily
$\gamma > s_+$, for otherwise $M_- (\alpha) \le B_\gamma \le M_+ (a) \le
W_-(a)+a$
which together with $\oo{b}(a)= -s_-(a)$ would imply that the depth of
$(\alpha, \oo{b} (a), \gamma)$ is smaller than $a$.  
Thus, under (\ref{sinai7}) if
$\alpha \le -s_-(a)$ then $\gamma > s_+$, and in this case
$\oo{b}(a) = s_+(a)$ since $B_{s_+(a)} < B_{-s_-(a)}$.  Finally, if
$\alpha> - s_-(a)$ then $\oo{b}(a) \not= -s_-(a)$ 
and hence $\oo{b}(a) = s_+(a)$.

Hence, we showed that if $W_-(a) > W_+(a)$ 
then (\ref{sinai7}) is equivalent to
$\oo{b}(a) = s_+(a)$.  Interchanging the positive and negative axis, we
conclude that if
$W_-(a) < W_+(a)$, then $\oo{b}(a)= -s_-(a)$ iff
$M_+(a) < (W_+(a) + 1) \vee M_+(a)$.  This completes the proof that 
$\oo{b}(a) = s_+(a)$ is equivalent to (\ref{sinai7}) or (\ref{sinai8}).

To complete the proof of the theorem, assume first
$W_-(a) > W_+(a)$.  Then, $\oo{b}(a) = s_+(a)$ iff (\ref{sinai7}) holds,
i.e.\ $M_+ (a) < (W_- (a) + a) \vee M_- (a) = H_-(a)$.
But $H_-(a) \ge W_- (a) + a \ge W_+ (a) + a$, and hence
$M_+(a) < H_-(a)$ is equivalent to $M_+(a) \vee (W_+(a) + a) < H_- (a)$,
i.e.\  $H_+(a) < H_- (a)$.
The case $W_+(a) < W_-(a)$ is handled similarly by using (\ref{sinai8}).
\qed

\section{Aging properties for Sinai's model}
One may use the representation in Theorem~\ref{the-sinai2} in order to
evaluate explicitly the law of $\oo{b}(a)$  (note that $\oo{b}(a)
\overset{{\cal L}}{=} a \oo{b} (1)$ by Brownian scaling).
This is done in \cite{kesten1}, 
and we do not repeat the construction here.
Our goal is to use Theorem~\ref{the-sinai2} to show that Sinai's model
exhibits \textsl{aging} properties.  More precisely, we claim that
\begin{theorem}
\label{the-sinai3}
Assume $P(\min(\om_0, \om_0^+) < \eps) = 0$ 
and Assumption~\ref{ass-sinai}.
Then, for $h>1$,
\begin{equation}
\label{sinai10}
\lim_{\eta\to0} \lim_{n\to\infty} \PP^o
\left(\frac{|X_{n^h} - X_n|}{(\log n)^2} < \eta\right)
= \frac{1}{h^2} \left[ \frac{5}{3} - \frac{2}{3} e^{-(h-1)} \right]
\end{equation}
\end{theorem}
Theorem \ref{the-sinai3} was derived heuristically
in \cite{fisher}, to which we refer for additional aging properties and 
discussion.  The right hand side of
formula \req{sinai10} appears also in \cite{golosov2}, in a slightly 
different context.

\proof
Applying Theorem~\ref{the-sinai}, the limit in the left hand side of
(\ref{sinai10}) equals
$$
\calQ\Bigl(\oo{b} (h) = \oo{b}(1)\Bigr) 
= 2\calQ \Bigl(\oo{b} (h) = \oo{b}(1) = s_+ (1) = s_+
(h)\Bigr)
$$
Note that
$$
\calQ(s_+ (h) = s_+(1)) = \calQ \left(
\begin{array}[t]{c}
\text{Brownian motion, started at height $1$,}\\
\text{hits $h$ before hitting $0$}
\end{array}
\right) = \frac{1}{h}\,.
$$
Hence, using that on $s_+(1) = s_+(h)$ one has
$W_+(1) = W_+(h), M_+(1) = M_+(h)$,
and using that the event $\{s_+(h)=s_+(1)\}$ depends only on increments 
of the path of the Brownian motion after time $T_+(1)$,
one gets
\begin{equation}
\label{sinai10a}
\calQ\Bigl(\oo{b}(h) = \oo{b} (1) \Bigr)
= \frac{2}{h} \calQ\Bigl( H_+ (1) < H_-(1), (W_+ (1) + h) \vee
M_+ (1) < H_-(h)\Bigr)\,.
\end{equation}
Next, let
\begin{align*}
\tau_0 & = \min\{t>s_-(1): B_{-t} = W_-(1) + 1\} \\
\tau_h & = \min \{t>\tau_0: \; B_{-t}= W_-(1) + h \,\text{\ or \ }
B_t = W_-(1)\}\,.
\end{align*}
Note that $\tau_h-\tau_0$ has the same law as that of the hitting time of
$\{0,h\}$ by a Brownian motion $Z_t$ with $Z_0=1$.  (Here, 
$Z_t = B_{-(\tau_0+t)} - W_-(1) !$).  Further, letting
$I_h = \won_{\{B_{\tau_h} = W_-(1)\}} 
(= \won_{\{Z_{\tau_h-\tau_0}=0\}})$,
it holds that
\begin{align*}
W_-(h) & = W_-(1) + I_h \tilde{W}_- (h) \\
M_-(h) & = \begin{cases}
M_-(1), & I_h=0\\
M_-(1) \vee (\oo{M}_- (h) + W_-(1) + 1) \vee (\tilde M_-(h) + W_-(1)),
& I_h=1
\end{cases}
\end{align*}
where $(\tilde{W}_-(h), \tilde{M}_-(h))$ are independent of 
$(W_-(1), M_-(1))$ and
possess the same law as $(W_-(h),$
$M_-(h))$, while $\oo{M}_-(h)$ 
is independent of
both $(W_- (1), M_-(1))$ and $(\tilde{W}_-(h), \tilde{M}_-(h))$ 
and has the law of the maximum of a Brownian motion, started 
at $0$, killed at hitting  $-1$ and conditioned not to hit  $h-1$.
(See figure 4 
for a graphical description of these random variables.)

\begin{figure}[h]
\label{sinai-fig4}
\begin{minipage}[b]{0.5\linewidth}
\centering
\includegraphics[scale=.5]{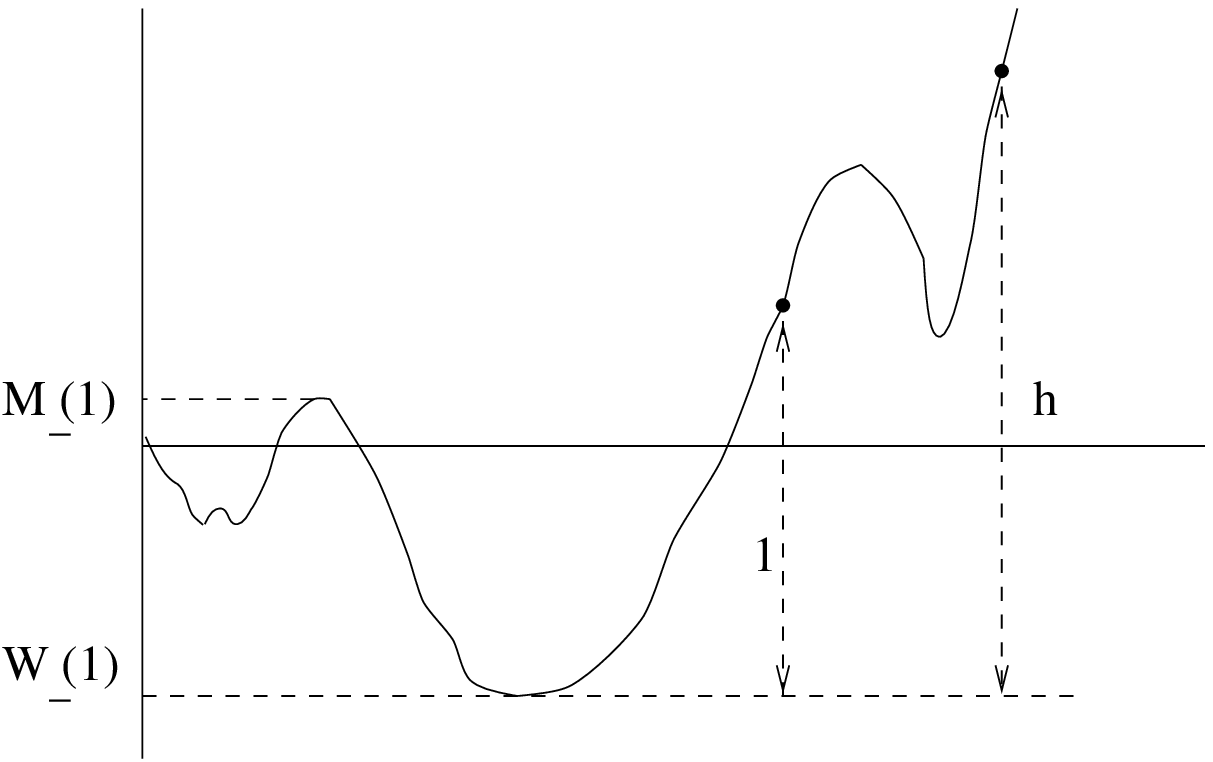}

{$I_h=0$} 
\end{minipage}
\begin{minipage}[b]{0.5\linewidth}
\centering
\includegraphics[scale=.5]{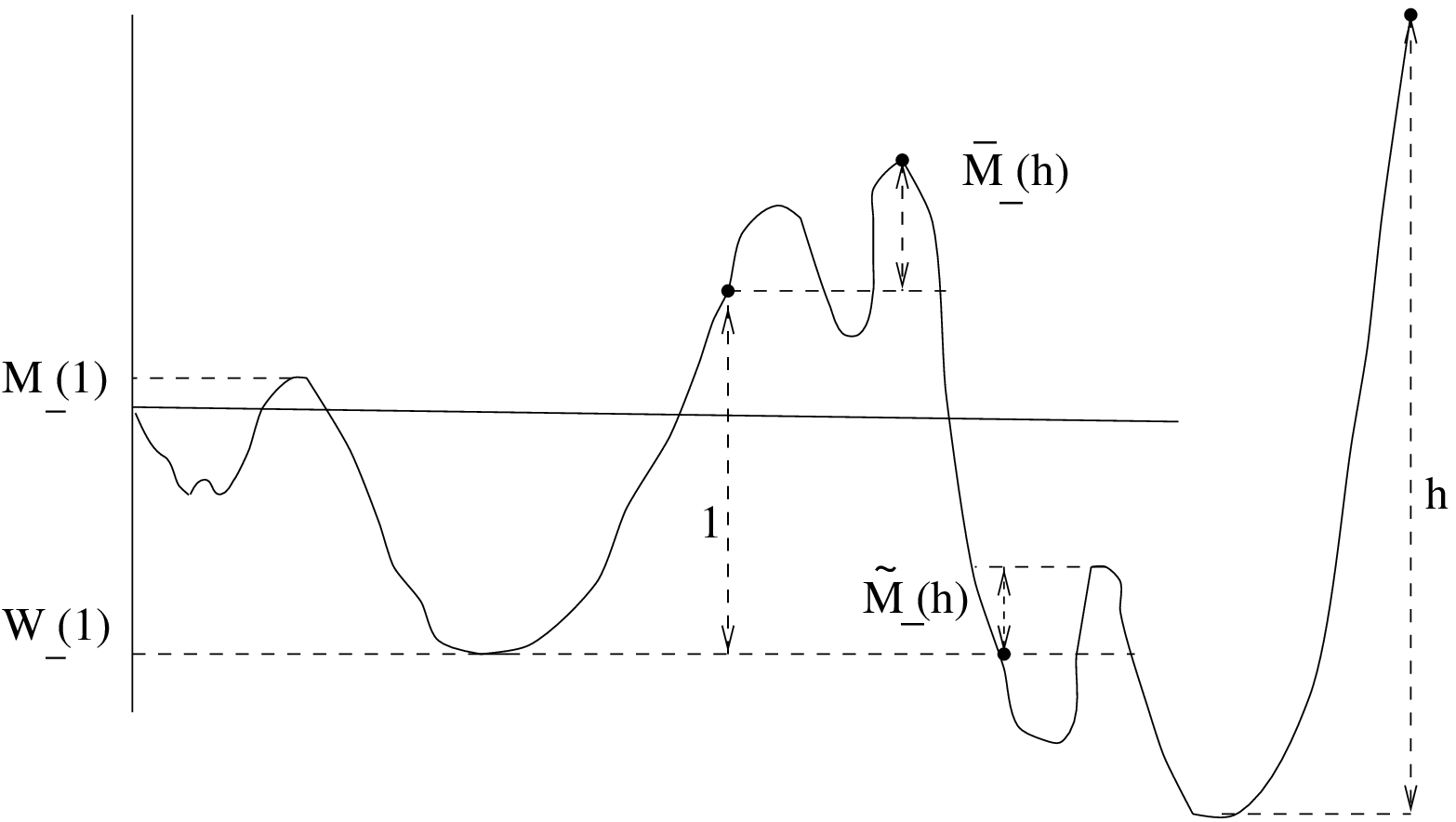}

{$I_h=1$}
\end{minipage}
\caption{Definition of auxiliary variables}
\end{figure}

Set now 
$$\hat{M}_-(h) = 
\begin{cases}
h, & I_h = 0\\
1+\oo{M}_-(h), & I_h=1\,,
\end{cases}
$$
$\tilde H_-(h)=(\tilde W_-(h)+h)\vee \tilde M_-(h)$
and $\Gamma(h) = \max (\tilde H_-(h), \hat M_-(h))$.  
Note that $\tilde H_-(h)$ has the same law as $H_-(h)$ but is 
independent of $\oo{M}_-(h)$.
Further, it is easy to check that
$(W_-(h)+h) \vee M_- (h) = (W_-(1) + \Gamma_h) \vee M_-(1)$
(note that either $M_-(h) = M_-(1)$ or
$M_-(h) > M_-(1)$ but in the latter case,
$M_-(h) \le W_-(1) + \Gamma(h)$.)
We have the following lemma, whose proof is deferred:
\begin{lemma}
\label{lem-sinai}
The law of
$\Gamma_h$ is $\frac{1}{h} \delta_h + \frac{h-1}{h} U [1,h]$, where
$U[1,h]$ denotes the uniform law on $[1,h]$.
\end{lemma}
Substituting in (\ref{sinai10a}), we get that
\begin{equation}
\label{sinai11}
\calQ
(\oo{b}(h) = \oo{b}(1)) = 
E_{\calQ}\left(E_{\calQ}(\oo{b}(h)=\oo{b}(1)|\Gamma(h))\right)=
\frac{2}{h^2} \left[ \int_1^h Q(t) dt + Q(h)\right]
\end{equation}
where
$$
Q(t) = \calQ(H_+(1) < H_-(1), H_+ (h) < H_-(t))\,.
$$
In order to evaluate the integral in (\ref{sinai11}), 
we need to evaluate the
joint law of $(H_+(1), H_+(t))$
(the joint law of $(H_-(1),H_-(t))$ being identical).
Since $0\le H_+(1) \le 1$ and $H_+(1) \le H_+(t) \le H_+(1) + t-1$,
the support of the law of $(H_+(1), H_+(t))$ is the domain A defined by
$0 \le x \le 1$, $x \le y \le x+t-1$.  Note that for
$(z,w) \in A$,
\begin{align*}
\calQ
(H_+(1) \le z, H_+(t) \le w) &=
\calQ(M_+(1) \le z \wedge w, W_+ (1) \le - [(1-z) \vee (t-w)]\\
& = \calQ\Bigl(M_+(1) \le z, W_+(1) \le - (t-w)\Bigr)
\,.
\end{align*}
We now have the following well know lemma. For completeness,
the proof is given at the end of this section:
\begin{lemma}
\label{lem-sinai2}
For $z+y \ge 1$, $0 \le z \le 1$, $y \ge 0$,
$$
\calQ(M_+(1) \le z, W_+(1) \le -y) = z e^{-(z+y-1)}\,.
$$
\end{lemma}
Lemma~\ref{lem-sinai2} implies that, for $(z,w) \in A$, $t>1$,
\begin{equation}
\label{sinai12}
\calQ(H_+(1) \le z, H_+(t) \le w) = z e^{-(z+t-w - 1)}\,.
\end{equation}
Denote by $B_1$ the segment $\{0\le x=y\le 1\}$ and by 
$B_2$ the segment
$\{t-1 \le y= x+t -1 \le t\}$.
We conclude, after some tedious computations, 
that the law of $(H_+(1), H_+(t))$:
\begin{itemize}
\item
possesses the density
$
f(z,\om) = (1-z) e^{-z} e^{-w-(t-1)},
\quad(z, w) \in A \backslash (B_1\cup B_2)
$
\item
possesses the density
$\tilde f(z,y) = (1-z) e^{-(t-1)}, \quad z=w \in B_1$
\item
possesses the density
$\oo{f} (z,z+t-1) = z, \quad w=z+t-1 \in B_2$.
\end{itemize}
Substituting in the expression for $Q(t)$, we find that
$$
Q(t) = \frac{5}{12} e^{-(h-t)} + \frac{1}{12} e^{-(h+t-2)}
\,.
$$
Substituting in (\ref{sinai12}), the theorem follows.
\qed

\noindent
{\bf Proof of Lemma \ref{lem-sinai}:}
Note that $\calQ(I_h=0)=1/h$, and in this case $\Gamma_h=h$.
Thus, we only need to consider the case where $I_h=1$ and show that
under this conditioning,
$\max(H_-(h), 1+\oo{M}_-(h))$ possesses the law $U[1,h]$.
Note that by standard properties of Brownian motion, 
$$\calQ(\hat M_-(h)\leq \xi|I_h=1)=\frac{\frac{\xi-1}{\xi}}{
\frac{h-1}{h}}.$$
We show below that the law  of $\tilde H_-(h)$, which is 
identical to the law 
of $H_-(h)$, is uniform on $[0,h]$. Thus, using independence,
for $\xi\in [1,h]$,
$$ \calQ(\Gamma_h<\xi|I_h=1)=
\frac{h(\xi-1)  \xi}{\xi (h-1) h}=
\frac{\xi-1}{h-1}\,,$$
i.e. the law of $\Gamma_h$ conditioned on $I_h=1$ is indeed $U[1,h]$.

It thus only remains to evaluate the law of $H_-(h)$. By 
Brownian scaling, the law of $H_-(h)$ is identical to the law of $hH_+(1)$,
so we only need to show that the law of $H_+(1)$ is uniform on $[0,1]$.
This in fact is a direct consequence of Lemma \ref{lem-sinai2}.
\qed

\noindent
{\bf Proof of Lemma \ref{lem-sinai2}: \ }
Let $\calQ^x$ denote the law of a Brownian motion
$\{Z_t\}$ starting at time $0$
at $x$. The Markov property now yields
\begin{eqnarray}
\label{sinaifinal}
\calQ(M_+(1) \le z, W_+(1) \le -y) &=& 
\calQ^o(\{Z_\cdot\} \,
\mbox{\rm hits $z-1$ before hitting $z$})
\calQ^{z-1}(M_+(1) \le z, W_+(1) \le -y)
\nonumber \\
&=&
z \calQ^o(M_+(1) \le 1, W_+(1) \le -y-z+1)\nonumber \\
& =& z \calQ^o(W_+(1)\leq -(y+z-1))\,.
\end{eqnarray}
For $x\geq 0$, let $f(x):=\calQ(W_+(1)\leq -x)$. The Markov
property now implies
$$f(x+\epsilon)=f(x) \calQ^{-x}(W_+(1)\leq -(x+\epsilon))=
f(x) f(\epsilon)\,.$$
Since $f(0)=1$ and $f(\epsilon)=1-\epsilon+o(\epsilon)$, 
it follows  that $f(x)=e^{-x}$. Substituting in 
\req{sinaifinal}, the lemma follows.
\qed



\end{document}